 \documentclass{amsart}
\usepackage{amsmath,amssymb,amsthm,verbatim,latexsym}
\usepackage[mathscr]{eucal}
\usepackage{amssymb}
\usepackage{graphicx}

\newcommand{\cB}{\mathcal{B}}
\newcommand{\cT}{\mathcal{T}}
\newcommand{\Pm}{\mathcal{P}}
\newcommand{\A}{\mathcal{A}}
\newcommand{\Z}{\mathbb{Z}}
\newcommand{\Yildirim}{Y{\i}ld{\i}r{\i}m}

\allowdisplaybreaks

\newtheorem{theorem}{Theorem}

\newtheorem{conjecture}{Conjecture C\hspace*{-1mm}}
\newtheorem*{hypothesisT}{Hypothesis T}

\theoremstyle{remark}
\newtheorem{remark}{Remark}

\newtheorem*{rema}{Remark}
\begin{document}

\title[Small gaps between almost primes]{S\lowercase{mall gaps between almost primes, the parity problem and some
conjectures of} E\lowercase{rd\H{o}s on consecutive integers\\
}}

\author{D. A. Goldston, S. W. Graham, J. Pintz and C. Y. Y{\i}ld{\i}r{\i}m}
\subjclass[2000]{Primary: 11N37; Secondary 11N36.}
\date{March 18, 2008}

\numberwithin{equation}{section}

\thanks
{The first and second authors were supported in part by NSF Grants,
the third author by OTKA grants No.  43623, 49693, 67676
and the Balaton program,
the fourth author by T\"UBITAK}

\begin{abstract}
In a previous paper, the authors proved that in any system of three linear forms
satisfying obvious necessary local conditions, there are at least two  forms
that infinitely often assume $E_2$-values; i.e., values that are products of exactly two primes. 
We use that result to prove that there are inifinitely many integers $x$ that 
simultaneously satisfy
$$
\omega(x)=\omega(x+1)=4, \Omega(x)=\Omega(x+1)=5, \text{ and } 
d(x)=d(x+1)=24.
$$
Here, $\omega(x), \Omega(x), d(x)$ represent the number of prime divisors of $x$,
the number of prime power divisors of $x$, and the number of  divisors of $x$,
respectively.
We also prove similar theorems where $x+1$ is replaced by $x+b$ for an arbitrary positive integer $b$. 
Our results sharpen earlier work of Heath-Brown, Pinner, and Schlage-Puchta.
\end{abstract}
\maketitle

\section{Introduction}

Erd\H{o}s had many favorite problems on consecutive integers (see
the work of Hildebrand \cite{Hil}). We will discuss among others the
celebrated Erd\H{o}s--Mirsky conjecture \cite{EM} on consecutive values
of the divisor function:

\begin{conjecture}
\label{C1}
$d(x)=d(x+1)$ infinitely often.
\end{conjecture}

We will also deal with the analogous conjectures for 
$\Omega(x)$ and $\omega (x)$, which denote the number of prime factors of
a positive integer $x$ counted with and without multiplicity,
respectively.

\begin{conjecture}
\label{C2}
$\Omega(x)=\Omega(x+1)$
infinitely often.
\end{conjecture}

\begin{conjecture}
\label{C3}
$\omega(x)=\omega(x+1)$
infinitely often.
\end{conjecture}

 The twin prime conjecture is
clearly equivalent to a stronger form of
Conjectures C\ref{C1}--C\ref{C2} for
$x+2$ in place of $x+1$ namely, with
\begin{equation} \label{eq:1.1}
d(x)=d(x+2)=2, \text{ and } \Omega(x)=\Omega(x+2)=1.
\end{equation}

However, the original
Conjectures C\ref{C1}--C\ref{C3} also follow (in a
stronger form, similarly to the above) from some analogues of the
twin prime conjecture. 
For example, it is conjectured that 
$2p+1$ is prime for infinitely many primes $p$,
and Chen's method
(see Chapter 11 of \cite{HR}) proves that there are infinitely many
$p\in\Pm$ ($\Pm$ denotes the set of prime numbers) such that either
\begin{equation}
2p+1 \in \Pm
\label{eq:1.2}
\end{equation}
or
\begin{equation}
2p+1={p_1}{p_2} \quad {p_1},{p_2}\in\Pm, \quad {p_1}\ne{p_2}.
\label{eq:1.3}
\end{equation}

If, as we believe,  \eqref{eq:1.3}  holds infinitely often, then
C\ref{C1}--C\ref{C3} hold; more precisely
\begin{equation}
d(2p)=d(2p+1)=4, \quad
\omega(2p)=\Omega(2p)=\omega(2p+1)=\Omega(2p + 1)=2.
\label{eq:1.4}
\end{equation}

Due to the connection to the  problem in \eqref{eq:1.3} (``which is
believed to be of the same depth as twin prime conjecture''
\cite[p.\ 309]{Hil}), Conjectures C\ref{C1}--C\ref{C3}
were also considered extremely
difficult, if not hopeless.

It was therefore a great surprise for Erd\H{o}s  \cite{Er} (and
probably for other number theorists as well) when C. Spiro \cite{Spi}
proved in 1981 that
\begin{equation}
d(x)=d(x+5040)\text{ infinitely often.}
\label{eq:1.5}
\end{equation}
At about the same time, Heath-Brown \cite{HB1} found a
conditional proof of C2 under a hypothesis slightly weaker than the
Elliott--Halberstam conjecture \cite{EH}. In 1984 he
succeeded \cite{HB2} in proving the original Erd\H{o}s--Mirsky
Conjecture C\ref{C1} by using the ideas of Spiro in combination of
other new ideas.
His method also yielded C2, but not C3.
About two decades later,  J.-C. Schlage--Puchta \cite{SchP} gave the first 
proof of C3.

A common feature of all these proofs was their intimate connection
with almost primes. More precisely (as suggested by
\eqref{eq:1.2}--\eqref{eq:1.3}),
all numbers produced that satisfied the relations
C\ref{C1}--C\ref{C3}  had a
bounded number of prime factors, but this number (even its parity) 
was left unspecified by the nature of the methods applied.

The phenomenon which ``prevents us from showing that \eqref{eq:1.3} has
infinitely many solutions'' is called the
``parity obstacle'', ``parity barrier" or ``parity problem".
Selberg gave examples of the two
sets \cite[Ch.\ 4]{Gr}
\begin{equation*}
{\A^{{(-1)}^r}}(X)=\{a: \; X<a\le2X, \; \Omega(a)\equiv r+1 \pmod{2}\}
\end{equation*}
which show that the upper and lower bounds obtained under general
conditions  of the linear sieve of Rosser are optimal. The
parity problem is expressed informally as saying that sieve methods
cannot differentiate  between integers with an even and an odd
number of prime factors. Therefore, the
parity obstacle prevents the sieve method from revealing the
existence of primes in a suitable set as formulated by Greaves \cite[p.\ 171]{Gr}.
For example, while Chen's method yields that $p+2$
(or $2p+1$) has infinitely often at most two prime factors, the
method is unable to specify the parity of the number of prime
factors of $p+2$ (or $2p+1$). As a result, the seemingly much
easier assertion that for infinitely many primes $p$, $p+2$ (or
$2p+1$) has an odd (or even) number of prime factors, is still
open \cite[p.\ 310]{Hil}.

The general view about the parity problem and
Conjectures C\ref{C1}--C\ref{C3} can be described  by again citing
the survey paper of Hildebrand \cite[p.\ 310]{Hil}:
``However, there is one crucial difference which makes this conjecture 
[the Erd\H{o}s--Mirsky conjecture] more accessible than 
`twin-prime type' conjectures. Namely, in contrast to the above-mentioned
problem on the parity of the number of prime factors of $2p + 1$ or
$p+2$ when $p$ is prime (and thus, in particular, has an odd
number of prime factors), when trying to prove that $d(n)=d(n+1)$
holds infinitely often one does not need to specify the parity of
the number of prime factors of $n$. It is this fact that allowed
the solution of the Erd\H{o}s--Mirsky conjecture, while bypassing
the deeper problems related to the twin prime conjecture.''

The situation is described in exactly the same way by Heath-Brown
in his work, where he gives the conditional solution of
C\ref{C2} \cite{HB1}:
``It should however be noted at this point that in solving 
$\Omega(n)=\Omega (n+1)$ we shall not have specified $\Omega (n)$, or
even the parity of $\Omega (n)$. Thus we avoid the parity problem,
rather than solve it.''
In his work \cite{HB2}, where he gives an unconditional solution of
C\ref{C1} and C\ref{C2}, he writes about the proof of
\eqref{eq:1.5} by Spiro (and his
remark refers also for the solution he gives for
C\ref{C1}--C\ref{C2}): 
``Thus one does not know the value of $\Omega (n)$ for the particular $n$
which satisfies $d(n)=d(n+5040)$.
In this way one sidesteps the `parity problem.' ''

In the present work we shall show that the method which yielded the
existence of short gaps between primes \cite{GPY1} and $E_2$-numbers (numbers with
exactly two distinct  prime factors) with bounded differences
\cite{GGPY1}, \cite{GGPY2} is able to show a stronger variant of
Conjectures C\ref{C1}--C\ref{C3}.
In this variant, the parity problem is not bypassed but overcome. We can show
 the above conjectures in the stronger form where  the value
of the relevant arithmetic function $d, \Omega$ or $\omega$ is specified.

\begin{theorem}
\label{th:1}
For any integer $A \ge 4$,  there are infinitely many integers $x$ with
\begin{equation}
\Omega(x)=\Omega(x+1)=A.
\label{eq:1.7}
\end{equation}
\end{theorem}

\begin{theorem}
\label{th:2}
For any integer $A \ge 3$,  there are infinitely many integers $x$ with
\begin{equation}
\omega(x)=\omega(x+1)=A.
\label{eq:1.8}
\end{equation}
\end{theorem}

\begin{theorem}
\label{th:3}
For any positive integer $A$ with $24 | A$,  there are infinitely many integers $x$ with
\begin{equation}
d(x)=d(x+1)=A.
\label{eq:1.9}
\end{equation}
\end{theorem}

As mentioned earlier, Heath-Brown \cite{HB2} was the first to prove Conjectures C1 and  C2, 
and Schlage-Puchta \cite{SchP} was the first to prove Conjecture C3. 
Here, we will give the first proof that there are infinitely many $x$ that simultaneously satisfy 
the conditions of Conjectures \ref{C1} through \ref{C3}.

For a positive integer $n$ with prime
factorization 
\begin{equation*}
n=p_1^{\alpha_1} p_2^{\alpha_2} \ldots p_k^{\alpha_k},
\end{equation*}
we define the {\it exponent pattern} of $n$ to  be the multiset 
$\{\alpha_1, \alpha_2, \ldots, \alpha_k\}$.  We follow the usual conventions for
multisets;  in particular,
an element may appear more than once, 
and the order of the elements does not matter.
This definition is relevant because if $x$ and $x+1$ 
have the same exponent pattern,
then $\omega(x)=\omega(x+1), \Omega(x)=\Omega(x+1),$ and $d(x)=d(x+1)$. 
We prove

\begin{theorem}
\label{th:123}
There exist infinitely many integers $x$ such that $x$ and $x+1$ both have exponent pattern
$\{2,1,1,1\}$. Consequently, there exist infinitely many integers $x$ such that
\begin{equation*}
\omega(x)=\omega(x+1)=4, \ \ \Omega(x)=\Omega(x+1)=5, \ \, 
\text{ and } d(x)=d(x+1)=24. 
\end{equation*}
\end{theorem}

We also note that in fact, we have $f(x)=f(x+1)$ infinitely often  for all 
 number-theoretic functions $f$ with the property  that $f(n)$ depends only on the exponent pattern 
 of $n$. Similar comments apply to Theorems  \ref{th:8a} and \ref{th:10a} below.

All of the theorems we prove here are done in a relatively straightforward and uniform 
way. In fact, all of our theorems are 
simple consequences of the
following result proved in \cite{GGPY2}.

\medskip

\noindent\textbf{Basic Theorem.} 
{\it  
We say that 
\begin{equation}
L_i(x)=a_i x+b_i \quad (1 \le i \le 3) \quad a_i,\;  b_i \in \Z,
\; a_i>0,
\label{eq:1.26}
\end{equation}
is an admissible triplet of linear forms if 
for every prime $p$ there exists $x_p \in \Z$ such that
\begin{equation*}
p \nmid L_1(x_p) L_2(x_p) L_3(x_p).
\end{equation*}
Let $C$ be any constant. 
If $\{L_1, L_2, L_3\} $ is an admissible triplet, then 
there are two forms $L_i(x), L_j(x)$ $(i\neq j)$ in the triplet
that simultaneously  take $E_2$-values with both prime factors
exceeding  $C$ for infinitely many integer values $x$.
In other words, there are two forms $L_i,L_j$ such that 
\begin{equation*}
\omega(L_i(x))=\Omega (L_i (x))= \omega (L_j(x))=\Omega (L_j(x))=2,
\end{equation*}
and 
\begin{equation*}
\Big( \prod_{p\le C} p, L_i(x) L_j (x) \Big)=1
\end{equation*}
for infinitely many integers $x$.
}

For future reference, we  make one further remark about the 
admissibility hypothesis in the Basic Theorem.
In order for the triplet $\{L_1,L_2, L_3\}$
to be admissible, it is obviously necessary to have
\begin{equation} \label{eq:1.27}
(a_i,b_i)=1, \ \ i=1,2,3.
\end{equation}
We shall refer to any system of linear forms that satisfies \eqref{eq:1.27} as a {\it reduced system.}
For any reduced system, 
\begin{equation*}
L_1(n) L_2(n) L_3(n) \equiv 0 \pmod p
\end{equation*}
has at most $3$ solutions $\pmod p$. 
Therefore, for any reduced system, 
the admissibility hypothesis is satisfied for any prime $p\ge 5$, and we need
check only  $p=2$ and $p=3$.

To illustrate how the Basic Theorem can be used, 
we give a short proof of Theorem \ref{th:2} in the case $A=3$. 
Consider the  system
\begin{equation} \label{eq:51}
L_1(m)=6 m+1, \; L_2(m)=8 m+1, \; L_3(m)=9m+1.
\end{equation}
This system is admissible because we may take $x_p=0$ for all primes $p$.
We note the relations
\begin{equation}\label{eq:52}
4L_1(m)=3L_2(m)+1, \ \  3L_1(m)=2L_3(m)+1, \ \ 9L_2(m)=8L_3(m)+1.
\end{equation}
By the Basic Theorem, at least two of the forms will be
simultaneously $E_2$-numbers for infinitely many values of $x$. 
We take $C=3$,
so that resulting $E_2$-numbers  have both prime factors exceeding $3$.

Now suppose the two forms giving infinitely many $E_2$-numbers 
are $L_1$ and $L_2$. Then we let
\begin{equation*}
x=3L_2(m), \ \  x+1=4L_1(m).
\end{equation*}
If the two relevant forms  are $L_1$ and $L_3$, we let
\begin{equation*}
x=2L_3(m), \ \  x+1=3L_1(m).
\end{equation*}
If the two relevant forms are  $L_2$ and $L_3$, we let
\begin{equation*} 
x=8L_3(m), \ \  x+1=9L_2(m).
\end{equation*}
In all cases, we obtain infinitely many positive integers $x$ with 
$\omega(x)=\omega(x+1)=3$.

As the proof illustrates, our approach is to combine the Basic Theorem with a judicious choice
of linear forms. 
Another element of our approach is the idea of  ``adjoining'' extra prime factors. 
We  will discuss this further in Section~\ref{sec:2};  e.g., see
the proof of Theorem~\ref{th:1}.

Next, we consider a modification of Conjectures \ref{C1} through \ref{C3} with $x, x+1$
replaced by $x, x+2$.  As noted before, this is directly related to the twin prime conjecture.
We emphasize that while we are unable to specify the parity of
$\omega (p+2)$ or $\Omega (p+2)$ for any infinite set of primes $p$, 
we can specify the common
value of $\omega (n)$ and $\omega (n+2)$ (or $\Omega (n)$ and
$\Omega (n+2)$), where the value can be given  almost
arbitrarily.

\begin{theorem}
\label{th:4}
For any $B \ge 4$ there are infinitely many integers $x$ with
\begin{equation}
\omega(x)=\omega(x+2)=B .
\end{equation}
\end{theorem}

\begin{theorem}
\label{th:5}
For any $B \ge 5$ there are infinitely many integers $x$ with
\begin{equation}
\Omega(x)=\Omega(x+2)=B .
\end{equation}
\end{theorem}

In the previous papers of this series \cite{GGPY1} \cite{GGPY2} we
investigated the distribution of $E_2$-numbers. We are unaware
of earlier results 
(or even earlier known methods capable of yielding such results) 
on $E_2$-numbers showing
\begin{equation}
\liminf_{n \to \infty} \frac{q_{n+1}-q_n}{\log q_n / \log \log q_n}
=0,
\label{eq:1.12}
\end{equation}
where $q_n$ denotes the $n^{\text th}$ $E_2$ number. 
This is in contrast to $\Pm_2$-numbers
(being the product of at most two primes) where Chen's method
proved that we may have infinitely often $n,\, n+2 \in \Pm_2$. 
The parity problem seems to indicate that problems with $E_2$-numbers
may have about the same difficulty as the analogous problems with
primes. We proved unconditionally that   
\begin{equation}
\liminf_{n \to \infty}(q_{n+1}-q_n)\le 6;
\label{eq:1.13}
\end{equation}
in fact, this follows from the Basic Theorem by considering the forms $\{n,n+2,n+6\}$. 
Although \eqref{eq:1.13} is clearly much stronger than
\eqref{eq:1.12}, it is still unclear whether there is any theoretical obstacle which
prevents us to find infinitely many pairs of $E_2$-numbers with a
fixed difference.

Our Theorems \ref{th:1}--\ref{th:5}
show, however, that
these barriers may be overcome if we increase (slightly) the
fixed number of prime divisors. Let  $r_n$ and $R_n$ 
be the $n^{\text{th}}$ positive integer with 
\begin{equation}
\omega (r_n)=3 \quad \text{and} \quad  \Omega(R_n)=3, \text{
respectively.}
\label{eq:1.14}
\end{equation}
An immediate corollary of Theorem~\ref{th:2} is that 
\begin{equation*}
\liminf_{n \to \infty} (r_{n+1}-r_n)=1.
\end{equation*}
We also prove 
\begin{theorem}
\label{th:6}
\begin{equation*}
\liminf_{n \to \infty} (R_{n+1}-R_n) \le 2.
\end{equation*}
\end{theorem}

Analogously to \eqref{eq:1.14}, let  $s_n$ and $S_n$ denote the 
$n^{\text{th}}$ positive integer  with
\begin{equation*}
\omega (s_n)=4 \quad \text{and} \quad  \Omega(S_n)=4, \text{ respectively.}
\end{equation*}
Theorems \ref{th:1} and \ref{th:2} imply
\begin{equation*}
\liminf_{n \to \infty} (s_{n+1}-s_n)= \liminf_{n \to \infty}
(S_{n+1}-S_n)=1.
\end{equation*}
The same holds if we consider the sequence ${s_n^{(\nu)}},
{S_n^{(\nu)}}$ of integers defined for an arbitrary integer $\nu \ge
4$ by
\begin{equation*}
\omega ({s_n^{(\nu)}})=\nu, \quad \Omega ({S_n^{(\nu)}})=\nu. 
\label{eq:1.17}
\end{equation*}

These results settle almost completely the problem of small
differences between consecutive almost primes having exactly $\nu \ge 3$ prime factors 
(except for the case when prime factors are
counted with multiplicity, and $\nu =3$ -- see Theorem~\ref{th:6}).

We also consider the generalization of  
Conjectures C\ref{C1}--C\ref{C3} to an arbitrary shift.
 In 1997,  Pinner \cite{Pn}
used  an ingenious extension of Heath-Brown's method
to prove that there are infinitely many positive integers $x$ with
\begin{equation*}
d(x)=d(x+n)
\end{equation*} 
for any fixed $n$. 
His method also works for $\Omega$, but it does not work for $\omega$.
However, Buttkewitz \cite{B} has recently proved that $\omega(x)=\omega(x+n)$
holds infinitely often for an infinite set $\cB$ of shifts $n$.

We will show that our method provides a positive answer to a
stronger form of these more general conjectures as well.
Similiarly to Theorems \ref{th:1}--\ref{th:5},
we can prescribe the value of the
$\omega$ and $\Omega$ functions in case of an arbitrary given shift $n$.
Our first results in this direction treat the case of even $n$. 

\begin{theorem}
\label{th:8a}
Let $n$ be an even positive integer. Then there exist infinitely many  $x$
such that $x$ and $x+n$ both have exponent pattern
$$ \{2,1,1,1,1\}.$$
In particular, there exist infinitely many $x$ such that
\begin{equation*}
\omega(x)=\omega(x+n)=5,\ \ \Omega(x)=\Omega(x+n)=6, \ \ 
\text{and } d(x)=d(x+n)=48.
\end{equation*}
\end{theorem}

\begin{theorem}
\label{th:9a} 
If $n\equiv 0 \pmod 4$, then for any integer $A\ge 4$, there
are infinitely many $x$ such that 
\begin{equation}\label{eq:1.21}
\Omega(x)=\Omega(x+n)=A.
\end{equation}
If $n\equiv 2 \pmod 4$, then for any integer $A\ge 5$, there are
infinitely many $x$ such that
\begin{equation}\label{eq:1.22}
\Omega(x)=\Omega(x+n)=A.
\end{equation}
If $n$ is even, then for any positive integer $A\ge 4$, there are infinitely many  $x$ such that
\begin{equation}\label{eq:1.23}
\omega(x)=\omega(x+n)=A.
\end{equation}
If $n$ is even, then for any  positive integer $A$ with $48|A$, there are infinitely
many $x$ such that 
\begin{equation}\label{eq:1.24}
d(x)=d(x+n)=A
\end{equation}
\end{theorem}

The proof of Theorem~\ref{th:8a} uses the linear forms
\begin{equation} \label{eq:1.25}
L_1(m)=2m+n, \ \ L_2(m)=3m+n, \ \ L_3(m)=5m+n.
\end{equation}
In an earlier draft of this paper, we used the same linear forms to show 
that if $n$ is odd and $n\not\equiv 15 \pmod {30}$, then there are infinitely
many $x$ such that $x$ and $x+n$ both have exponent pattern 
$\{2,1,1,1,1\}$. We eventually realized that one critical property of the forms
in \eqref{eq:1.25} is that $(3,5)$ is a twin prime pair. We then discovered 
Theorem \ref{th:10a}, which applies to any $n$ satisfying the following 

\begin{hypothesisT} Assume that $n$ is an odd positive integer and 
that there exists a twin prime pair $(p, p+2)$ such that 
$p(p+2)\nmid n$. 
\end{hypothesisT}

\begin{theorem} \label{th:10a}
If $n$ satisfies Hypothesis T, then there exist infinitely many positive integers $x$ 
such that $x$ and $x+n$ both have exponent pattern 
$$ \{2,1,1,1,1\}.$$
In particular, there are infinitely many $x$ such that 
\begin{equation*}
\omega(x)=\omega(x+n)=5, \ \ \Omega(x)=\Omega(x+n)=6, \ \ \text{and } 
d(x)=d(x+n)=48.
\end{equation*}
\end{theorem}

By combining Theorem~\ref{th:10a} with the trick of adjoining
prime factors, we get an analog of Theorem~\ref{th:9a} that holds 
for $n$ satisfying Hypothesis T.

\begin{theorem} \label{th:11a}
Assume that $n$ satisfies Hypothesis T.
For any positive integer $A\ge 5$, there are infinitely many $x$ such that
\begin{equation} \label{eq:11o}
 \omega(x)=\omega(x+n)=A.
 \end{equation}
For any positive integer $A$ with $48|A$, there are infinitely many $x$ such that
\begin{equation} \label{eq:11d}
d(x)=d(x+n).
\end{equation}
\end{theorem}

Along the same lines, we can
also prove that $ \Omega(x)=\Omega(x+n) = A$
for infinitely many $x$, any $n$ satisfying Hypothesis T, and any $A\ge 6$.
However, we will prove a stronger result on $\Omega$ in Theorem \ref{th:12a}.

The previous two theorems lead to the natural problem of determining 
which integers satisfy Hypothesis T.
Take any set $\cT_1$ with the property that 
$p\in \cT_1$ implies that $(p,p+2)$ is a twin prime pair. Let $\cT_2=\{p+2:p\in \cT_1\}$,
$\cT=\cT_1\cup \cT_2$, and 
\begin{equation} \label{eq:S}
 S=S_{\cT} =  \prod_{p\in \cT} p.
 \end{equation}
Then any $n\not\equiv S \pmod {2S}$ satisfies Hypothesis T. 

For example, if we take $\cT_1=\{3,5\}$, we get $S=105$, and   so any odd  $n$ with 
$n\not\equiv 105 \pmod{210}$ satisfies Hypothesis T.
Of course, one may take larger sets as well.%
\footnote{Note that $\cT_1$ and $\cT_2$ are disjoint if  $3\notin \cT_1$ and $5\notin \cT_1$.} 
Let 
\begin{equation} \label{eq:T1}
\cT_1 =\{ p: 10^{15} < p \le 10^{16} , (p, p+2)  \text{ is a twin prime pair}  \},  
\text{ and }  
 S =\prod_{p\in\cT_1} p(p+2).
\end{equation}
According to calculations of Sebah \cite{Seb} \cite{W}, the set $\cT_1$ 
has more than $9 \cdot 10^{12}$ elements,
and therefore 
$$ \log_{10} S > 30 \cdot 9\cdot 10^{12} >2\cdot 10^{14}.$$
Consequently, Hypothesis T is true for any odd $n$ with $n<10^{2\cdot 10^{14}}.$
Moreover, the density of $n$ which fail Hypothesis T is $\le 1/(2S)$, which is 
extremely small. 

In fact, Hypothesis T is probably true for all odd $n$, but this appears difficult to prove.
The truth of Hypothesis T for all odd $n$  is equivalent 
to the existence arbitrarily large sets $\cT$, and this in turn is equivalent to the 
twin prime conjecture. In light of this, it is desirable to give unconditional proofs 
of  \eqref{eq:11o} and  \eqref{eq:11d} for all odd $n$. 

\begin{theorem} \label{th:12a}
Assume that $n$ is an odd positive integer.
For any integer $A\ge 5$, there are infinitely many $x$ such that
\begin{equation} \label{eq:12O}
\Omega(x)=\Omega(x+n) = A.
\end{equation}
For any integer $A\ge 6$, there are infinitely many $x$ such that 
\begin{equation} \label{eq:12o}
 \omega(x)=\omega(x+n)=A.
 \end{equation}
For any positive integer $A$ with $288|A$, there are infinitely many $x$ such that
\begin{equation} \label{eq:12d}
d(x)=d(x+n)=A.
\end{equation}
\end{theorem}

In the  proofs of \eqref{eq:12o} and \eqref{eq:12d}, we employ the system 
$$ L_1(m)= 672m+41,\ \ L_2(m)=672m+47, \ \ L_3(m)=672m+55.$$
The Basic Theorem says that at least two of these forms represent $E_2$-numbers
infinitely often. This provides a usable analogue of the twin prime conjecture.
The proof of \eqref{eq:12O} is simpler; it uses the linear forms given in \eqref{eq:1.25}.

\section{The case of consecutive integers--Theorems \ref{th:1} through \ref{th:6}}
\label{sec:2}

For the proof of Theorem \ref{th:1}, we begin with the admissible system
\begin{equation}
L_1(m)=4m+1,\ \ L_2(m)=5m+1,\ \ L_3(m)=6m+1,
\label{eq:54}
\end{equation}
and we note  the relations
\begin{equation}
5L_1 =  4L_2+1, \ \  3L_1=2L_3+1, \ \ 6L_2=  5L_3+1.
\label{eq:55}
\end{equation}

If we apply the Basic Theorem directly to this system, we fail to get the desired result. 
If, for example, $L_1(m)$ and $L_2(m)$ are both $E_2$ numbers, and if $x=4L_2(m)$, 
then $\Omega(4L_2(m))=4$ and $\Omega(5L_1(m))=3$. 
To remedy this, we modify our system so that $4m+1$ and $6m+1$ each have one additional
prime factor.  We use the Chinese Remainder Theorem to solve the system of congruences
\begin{equation*}
4m+1  \equiv 0 \pmod{11}, \ \ 
6m+1  \equiv 0 \pmod{7}.
\end{equation*}
The solution is $m\equiv 8 \pmod {77}$.
Accordingly, we take
\begin{equation}
K(\ell)= 77\ell+8, \label{eq:57}
\end{equation}
and we define three new forms $K_1,K_2,K_3$ by setting $r_1=11, r_2=1, r_3=7$, and  
\begin{equation*}
K_i(\ell)= \frac{L_i(K(\ell))}{r_i} \text{\quad} (i=1,2,3).
\end{equation*}
In other words, 
\begin{equation} \label{eq:59}
K_1(\ell)=28\ell+3, \ \ K_2(\ell)=385\ell+41, \ \ K_3(\ell)=66\ell+7.
\end{equation}
This is an admissible system; we may use the comment after
\eqref{eq:1.27} for $p>3$. Otherwise, we  let $x_2=0$ and $x_3=2$. 
Therefore, we have two forms in \eqref{eq:59} that
are $E_2$-numbers with all prime factors exceeding $11$.
Consequently, there are at
least two relations among
\begin{equation*}
\Omega (4m+1)=3, \ \  \Omega(5m+1)=2, \ \  \Omega (6m+1)=3,  
\end{equation*}
that are simultaneously true for infinitely many integer values
$m$. In view of the relations \eqref{eq:55} we obtain
infinitely many integers $x$ with
\begin{equation}
\Omega (x)= \Omega (x+1)=4, \label{eq:61}
\end{equation}
since the multiplications in \eqref{eq:55} adjoin one new prime factor to
both $4m+1$ and $6m+1$  and two prime factors to $5m+1$. 

For the case $A>4$, we modify the above procedure slightly.
We let 
$$r_1= 11\cdot 13^{A-4}, \ \  r_2=17^{A-4}, \ \ r_3=7\cdot 19^{A-4}.$$
and we use the Chinese Remainder Theorem to find  $k $ such that 
\begin{equation*}
4k+1  \equiv  0 \pmod{r_1}, \ \ 
5k+1  \equiv  0 \pmod{r_2}, \text{ and } 
6k+1  \equiv  0 \pmod{r_3}.
\end{equation*}
We then take $K(\ell)= r_1 r_2 r_3 \ell+ k$, and define the forms $K_1, K_2, K_3$
by 
\begin{equation*}
K_i(\ell)= \frac{L_i(K(\ell))}{r_i} \text{\quad} (i=1,2,3).
\end{equation*}
We apply the Basic Theorem with $C=19$ to the system 
$\{K_1, K_2, K_3\},$ and the desired result follows.

\begin{rema}
We will often apply the above procedure of adjoining new prime factors
in the following, but we will usually omit the details.
\end{rema}

We have already proved Theorem \ref{th:2} in the case $A=3$; see the discussion beginning 
with equation \eqref{eq:51}. 
We do the case $A>3$ by the procedure of adjoining new prime 
factors. We take $r_1, r_2, r_3$ to be pairwise relatively prime
squarefree numbers with exactly $A-3$ prime factors all exceeding 3, and we start with the 
system \eqref{eq:51} as in the case $A=3$. 

We will derive Theorem \ref{th:3} from Theorem \ref{th:123}, so we now turn our attention to the latter. 
We begin with the linear forms
\begin{equation*}
L_1(m)=3m+2, \ \ L_2(m)=4m+3, \ \ L_3(m)=10m+7.
\end{equation*}
We note the relations
\begin{equation*}
3L_2=4L_1+1, \ \ 3L_3=10L_1+1, \ \ 5L_2=2L_3+1,
\end{equation*}
and we take 
\begin{equation*}
r_1=5,\ \  r_2=7^2, \ \ r_3 = 11^2. 
\end{equation*}
We see that 
$$L_i(m)\equiv 0 \pmod{r_i}  \text{\quad} (i=1,2,3) $$
holds for  $m\equiv 3956 \pmod {29645}$. 
Let $K(m)=29645m+3956$, and consider the three reduced forms
$$K_i(m) = \frac{ L_i(K(m))}{r_i}.$$
A calculation reveals that
\begin{equation*}
K_1(m)=17787m+2374, \ \ K_2(m)=2420m+323, \ \ K_3(m)=2450m+327.
\end{equation*}
The above three forms give an admissible system. To see this, we
use the comment  after \eqref{eq:1.27}
for primes $p\ge 5$. For $p=2$ and $p=3$, we note that 
\begin{equation*}
5\cdot 7^2\cdot 11^2\cdot K_1K_2K_3 (1) = L_1L_2L_3(K(1)) \equiv L_1L_2L_3(1) \equiv 1 \pmod 6,
\end{equation*}
so we may take $x_2=x_3=1$.
We apply the Basic Theorem and deduce that at least two of these forms 
are infinitely often simultaneously $E_2$-numbers with all prime factors 
exceeding 11. If $K_1(m)=p_1p_2$ and $K_2(m)=p_3p_4$, then we may take
\begin{equation*}
x=4L_1(m)=2^2\cdot 5\cdot p_1p_2, \ \ x+1=3\cdot L_2(m)=3\cdot 7^2\cdot p_3p_4.
\end{equation*}
The other cases are similar. We conclude that there are infinitely many $x$ such that
both $x$ and $x+1$ have exponent pattern
\begin{equation*}
\{2,1,1,1\}.
\end{equation*}
This proves Theorem \ref{th:123} and Theorem \ref{th:3} in the case $A=24$. 
When $A=24B$ with $B>1$, we modify the above with the procedure
of adjoining new prime factors. We can do this, for example, by using
\begin{equation*}
r_1=5\cdot 13^{B-1}, \ \ r_2=7^2\cdot 17^{B-1}, \ \ r_3= 11^2\cdot 19^{B-1}.
\end{equation*}
Note that this produces infinitely many $x$ such that both $x$ and $x+1$ have 
exponent pattern 
$$ \{2,1,1,1,B-1\}.$$
Consequently, we can say somewhat more; namely, there are infinitely many $x$ with
\begin{equation*}
d(x)=d(x+1)=24B, \ \  \Omega(x)=\Omega(x+1)=B+4, \ \ \text{and } 
\omega(x)=\omega(x+1)=5.
\end{equation*}

Theorems \ref{th:4} and \ref{th:5} are special cases of Theorem \ref{th:9a}, 
but we give the following independent short proofs.
For Theorem \ref{th:4} in the case $B=4$ and Theorem \ref{th:5} in the case $B=5$, consider the 
admissible system
\begin{equation*}
L_1(m)=2m+1, \ \ L_2(m)=12m+5, \ \ L_3(m)=20m+9.
\end{equation*}
Note the relations
\begin{equation*}
12L_1-2L_2=20L_1-2L_3=3L_3-5L_2=2.
\end{equation*}
Using the Basic Theorem and adjoining one new prime factor to each of $L_2$ and $L_3$,
 we get Theorem~\ref{th:4}  for $B=4$. Adjoining two new prime factors to each of $L_2$ and $L_3$
 gives Theorem~\ref{th:5} for $B=5$. The general cases may be 
 done by adjoining further prime factors.
 
 Finally, in order to prove Theorem~\ref{th:6}, consider the
admissible system
\begin{equation}
24 m + 1, \quad  36m + 1, \quad 72m + 1.
\label{eq:204}
\end{equation}
Using the Basic Theorem and adjoining one new prime factor to the
last form we obtain infinitely many positive integers $m$ such
that at least two of the relations
\begin{equation}
\Omega(24m + 1) = 2, \quad
\Omega(36m + 1) = 2, \quad
\Omega(72m + 1) = 3
\label{eq:205}
\end{equation}
hold. Equivalently at least two of the relations
\begin{equation}
\Omega(72m + 3) = 3, \quad
\Omega(72m + 2) = 3, \quad
\Omega(72m + 1) = 3
\label{eq:206}
\end{equation}
will be true for infinitely many positive integers $m$,
and this proves Theorem~\ref{th:6}.

\section{The case of even shift--Theorems \ref{th:8a} and  \ref{th:9a}}
\label{sec:3}

Throughout this section, we assume that $n$ is even, and we write $n=2n_2$. 

We consider the
three (non-admissible) forms
\begin{equation}
L_1(m)=2m+n, \ \  L_2(m)=3m+n, \ \ L_3(m)=5m+2n, \label{eq:141}
\end{equation}
with the relations
\begin{equation}
3L_1=2L_2+n, \ \  5L_1=2L_3+n,\ \   3L_3=5L_2+n.\label{eq:141a}
\end{equation}

\noindent We will restrict $m$ to some residue class $a\pmod {60}$, where
$a$ depends on $n$, 
\begin{equation} \label{eq:141c}
a \equiv   a_4 \pmod 4,  \ \ a \equiv  a_3 \pmod 3,  \ \ a\equiv  a_5 \pmod 5,
\end{equation}
and $a_4,a_3,a_5$ are defined as follows.  
First, we take
\begin{equation*}
a_4=
 \begin{cases}
   1 & \text{ if $n_2$ is even,}\\
   n_2 & \text{ if $n_2$ is odd.}\\
 \end{cases}
 \end{equation*}
 In other words, $a_4 \equiv 3 \pmod 4$ if $n_2\equiv 3\pmod 4$ and $a_4\equiv 1 \pmod 4$ 
 otherwise.
 If $3\nmid n$, we take $a_3=0$. If $3|n$,  we write $n=3n_3$ and choose $a_3$ so that
 $a_3(a_3+n_3) \not\equiv 0\pmod 3$. For example, we may take
\begin{equation*}
 a_3=
 \begin{cases}
        1 & \text{ if $n_3 \not\equiv 2 \pmod 3$,} \\
        2 & \text{ if $n_3 \equiv 2 \pmod 3$.}\\
 \end{cases}
\end{equation*}
 If $5\nmid n$, we take $a_5=0$. If $5|n$,  we write $n=5n_5$ and choose 
 $a_5$ so that $a_5(a_5+2n_5)\not\equiv 0 \pmod 5$. For example, we may take
 \begin{equation*}
 a_5=
 \begin{cases}
        1 & \text{ if $n_5 \not\equiv 2 \pmod 5$,} \\
        2 & \text{ if $n_5 \equiv 2 \pmod 5$.}\\
 \end{cases}
\end{equation*}

With this choice of $a$ we can show that if $m\equiv a \pmod {60}$, then 
 \begin{align}
 2^1 \parallel L_1(m)  \text{ if $4|n$,\quad}  &
2^2 \parallel L_1(m)  \text{ if $4\nmid n$,} \label{eq:142a}\\
 3^1 \parallel L_2(m)  \text{ if $3|n$, \quad}   &
3\nmid L_2(m)  \text{ if $3\nmid n$,} \label{eq:142b}\\
5^1 \parallel L_3(m)  \text{ if $5|n$,\quad} &
5\nmid L_3(m)  \text{ if $5\nmid n$.}\label{eq:142c}
 \end{align}
 To prove \eqref{eq:142a}, note that if $4|n$, then $2|n_2, a_4=1$ and 
$$ \frac 12 L_1(m) = m+n_2 \equiv a_4 \equiv 1 \pmod 2.$$
If $4\nmid n$, then  $n_2$ is odd, and  
$$\frac 12 L_1(m) = m+n_2 \equiv a_4+n_2 \equiv 2n_2 \equiv 2 \pmod 4.$$
To prove \eqref{eq:142b}, assume first that $3|n$. Then $n=3n_3$, and 
$$ \frac 13 L_2(m) = m+n_3 \not\equiv 0  \pmod 3.$$
If $3\nmid n$, then  $m\equiv 0 \pmod 3$ and 
$$ L_2(m) =  3m+ n \equiv n \not\equiv 0 \! \pmod 3.$$
The proof of \eqref{eq:142c} is similar to
\eqref{eq:142b}, and we leave the details to the reader.

Now let $p_i$ be the $i^{\rm th}$ prime, so that $p_1=2,p_2=3,p_3=5$. 
Let $b_1=b_2=1, b_3=2$. With this notation,
\begin{equation*}
L_i(m)= p_im+b_i n
\end{equation*}
We claim that if $m\equiv a \pmod {60}$, then 
\begin{equation} \label{eq:143}
p_i \nmid L_j(m)
\end{equation}
whenever $i\ne j$. 
For $i=2$ or $3$, this follows because $p_i$ divides either $m$ or $n$ but not both.
For $p_1=2$, this follows because $m$ is odd and $n$ is even.

We rephrase \eqref{eq:142a} through \eqref{eq:142c} by saying that  
if $m\equiv a \pmod {60}$, then 
\begin{equation} \label{eq:144}
2^{\alpha_1} \parallel L_1(m), \ \ 
3^{\alpha_2}\parallel L_2(m), \ \ 
5^{\alpha_3} \parallel L_3(m),
\end{equation}
where 
\begin{equation} \label{eq:145}
\alpha_1 = 
 \begin{cases}
  1 & \text{ if $4|n$},\\
  2 & \text{  if $4\nmid n$},\\
\end{cases} 
\ \ 
\alpha_2 = 
  \begin{cases}
  1 & \text{ if $3|n$},\\
  0 & \text{  if $3\nmid n$},\\
 \end{cases}
 \ \ 
\alpha_3 = 
  \begin{cases}
  1 & \text{ if $5|n$},\\
  0 & \text{  if $5\nmid n$}.\\
 \end{cases}
\end{equation}

We adjoin extra factors to the $L_i$ according to the following recipe. 

If $\begin{cases} 4|n \\ 4\nmid n \end{cases} $, then we adjoin a factor of 
 $\begin{cases} 19^2 \\ 19  \end{cases}$ to $L_1$.

If $\begin{cases} 3|n \\ 3\nmid n \end{cases} $, then we adjoin a factor of 
 $\begin{cases} 7^2 \\ 7^2\cdot 13  \end{cases}$ to $L_2$.
 
 If $\begin{cases} 5|n \\ 5\nmid n \end{cases} $, then we adjoin a factor of 
 $\begin{cases} 11^2 \\ 11^2\cdot 17  \end{cases}$ to $L_2$.

By doing this, we have arranged so each $L_i$ has a fixed divisor
of the form $r_i q_i^2$, where
\begin{equation*}
r_1=\begin{cases} 2 & \text{ if $4|n$} \\ 19 & \text{ if $4\nmid n$} \end{cases},
r_2=\begin{cases} 3 & \text{ if $3|n$} \\ 13 & \text{ if $3\nmid n$} \end{cases},
r_3=\begin{cases} 5 & \text{ if $5|n$} \\ 17 & \text{ if $5\nmid n$} \end{cases},
\end{equation*}
and 
\begin{equation*}
q_1=\begin{cases} 19 & \text{ if $4|n$} \\ 2 & \text{ if $4\nmid n$} \end{cases}, \,\,
q_2=7, q_3=11.
\end{equation*}
We may therefore apply the Basic Theorem to the three forms
\begin{equation} \label{eq:215}
K_i(\ell)  = \frac{L_i(60\ell+a)}{r_i q_i^2}.
\end{equation}
We obtain at least two different indices $i \in \{1,2,3\}$
such that the numbers in \eqref{eq:215} would be $E_2$-numbers
with both prime divisors  exceeding  $19$. Therefore, the numbers
$p_j (p_i m + b_i n)$ $(j \neq i)$  have the exponent-pattern
\begin{equation*}
\{2,\ 1,\ 1,\ 1,\ 1\},
\end{equation*}
This completes the proof of Theorem \ref{th:8a}.

Now we consider Theorem \ref{th:9a}. 
Statement \eqref{eq:1.24} follows from the construction in 
Theorem \ref{th:8a}; let $A=48K$ and adjoin 
a factor of $23^{K-1}$ to the forms $L_1, L_2, L_3$. 
For the proofs of the other statements, we take $L_1, L_2, L_3$ to
be  as in \eqref{eq:141}. Then we follow the proof of 
Theorem \ref{th:8a} down to \eqref{eq:145}.
Now let $K(\ell)=60\ell+a$, and set
\begin{equation*}
K_i(\ell) =\frac{ L_i(K(\ell))}{p_i^{\alpha_i}}
\end{equation*}
for $i=1,2,3$. The system $K_1,K_2,K_3$ is admissible, 
so we may apply the Basic Theorem.
By adjoining further primes as necessary, we see that there 
are infinitely many $x$ and at least two forms with
\begin{equation*}
\omega(p_j x+ b_j n) = 3, \ \ \Omega(p_j x+ b_j n) = 2+ \alpha_1,
\end{equation*}
where $\alpha_1 $ was defined in \eqref{eq:145}.
If we multiply any fixed $p_j x+ b_j n$ by $p_i  \; (i\ne j )$, then we obtain
two forms with 
\begin{equation*}
\omega(p_i(p_jx+b_jn)) = 4, \ \ \Omega(p_i(p_jx+b_jn)) = 3+\alpha_1.
\end{equation*}
This proves \eqref{eq:1.21} with $A=4$,  \eqref{eq:1.22} with $A=5$, and
\eqref{eq:1.23} with $A=4$. For larger values of $A$, we adjoin
further prime factors as necessary.

\section{Hypothesis T: Proofs of Theorems \ref{th:10a} and \ref{th:11a}} \label{sec:4}
In this section, we assume that $n$ satisfies Hypothesis T.
Consequently, $n$ is odd, and there is a twin prime pair $(p,p+2)$ such that
$p\nmid n$ or $(p+2)\nmid n$.
We divide our argument into two cases depending on whether or not $p=3$. 

{\bf Case 1.} $p=3$. In this case, we consider the three (non-admissible) linear forms
\begin{equation} \label{eq:400}
L_1(n) = 2m+n, \ \  L_2(m)=3m+n, \ \  L_3(m)=5m+2n.
\end{equation}
We note the relations 
\begin{equation} \label{eq:401}
3L_1-2L_2= 5L_1-2L_3=3L_3-5L_2=n.
\end{equation}
Our argument for this case is broken into three subcases: 
(a) $3\nmid n \text{ and } 5\nmid n$, 
(b) $3\nmid n \text{ and }  5|n$, 
(c) $3|n \text{ and } 5\nmid n$. 

{\bf Subcase 1a}. Assume that $3\nmid n$ and $5\nmid n$. 
We choose $m \pmod{180}$ by the relations
\begin{align}
m\equiv & 0 \pmod 4, \label{eq:402a}\\
m\equiv & 6+5n  \pmod 9, \label{eq:402b}\\
m\equiv & 0 \pmod 5. \label{eq:402c}
\end{align}
From the hypotheses,  $(n,30)=1$, so
$$ 2\nmid L_1(m), \ \ 3\nmid L_2(m), \ \ 5\nmid L_3(m).$$
Condition \eqref{eq:402a} implies that $m$ is even, so $ 2\nmid L_2(m)$.
However, $n$ is odd, so
$$ L_3(m)= 5m+2n \equiv 2 \pmod 4,$$
and  $2^1|| L_3(m)$. 
Condition \eqref{eq:402b} implies that 
$$L_1(m)= 2m+n \equiv 2n \not\equiv 0 \pmod 3
\text{ and }  L_3(m) = 5m+2n \equiv 3 \pmod 9,$$
so $3\nmid L_1(m) $ and $3^1|| L_3(m)$.
Finally, condition \eqref{eq:402c} implies that 
$$ L_1(m)  \equiv L_2(m)  \equiv n  \not\equiv  0 \pmod 5.$$
We deduce that  $L_1,L_2$ have no fixed divisors,  and
$L_3$ has the fixed divisors $2$ and $3$ with 
exponent one.
We adjoin  extra factors of $11^2$ to $L_1$ and  $13^2$ to $L_2$.
We apply the Basic Theorem to the reduced forms
\begin{equation*}
\frac{L_1(m)}{11^2}, 
\frac{L_2(m)}{13^2},
\frac{L_3(m)}{2\cdot 3}.
\end{equation*}
Multiplying the original forms by $2$, $3$ and $5$ according to
\eqref{eq:401},  
we arrive in all cases at infinitely many pairs  $x, x+n$ with exponent
pattern $\{2,1,1,1\}$. 
We adjoin one further prime factor to each
form to reach the exponent pattern $\{2,1,1,1,1\}$ enunciated in 
Theorem \ref{th:10a}.

{\bf Subcase 1b}. Assume that $3\nmid n$ and $5|n$.
Write $n=5n_5$. 
In this case,  we choose $m \pmod{180}$ by the relations
\begin{align}
m\equiv & 0~\pmod 4, \label{eq:403a}\\
m\equiv & 6+5n \pmod 9, \label{eq:403b}\\
m\equiv &
 \begin{cases}
      1 \pmod 5 & \text{ if $n_5\not\equiv 2 \pmod 5$},  \\
      2 \pmod 5 & \text{ if $n_5\equiv 2  \pmod 5$}. \\
  \end{cases}
  \label{eq:403c}
\end{align}
As in Subcase 1a, we have 
$$2\nmid L_1L_2(m), \ \ 3\nmid L_1L_2(m), \  \ 2^1||L_3(m), \ \ 3^1||L_3(m).$$
By \eqref{eq:403c}, 
$m(m+2n_5)\not\equiv 0 \pmod 5$, thus
$ 5\nmid L_1L_2(m) \text { and } 5^1|| L_3(m).$
We deduce that  $L_1,L_2$ have no fixed divisors  and
$L_3$ has the fixed divisors $2,3,5$ with 
exponent one.
We adjoin factors of $11^2 \cdot 17$ to $L_1$ and  $13^2 \cdot  19$ to $L_2$, then
we apply the Basic Theorem to the reduced forms
\begin{equation*}
\frac{L_1(m)}{11^2 \cdot 17}, \ \ 
\frac{L_2(m)}{13^2 \cdot 19}, \ \ 
\frac{L_3(m)}{2\cdot 3 \cdot 5}.
\end{equation*}
Multiplying the original forms by $2$, $3$ and $5$ according to
\eqref{eq:401},  we arrive in all cases at numbers with exponent pattern
$\{2,1,1,1,1\}$.

{\bf Subcase 1c.}  Assume that  $3|n$ and $5\nmid n$. 
Write $n=3n_3$.  This case is  similar to Subcase 1b; the primary difference is that 
the roles of $3$ and $5$ have been reversed.
In this case,  we choose $m \pmod{300}$ by the relations
\begin{align}
m\equiv & 2+ n \pmod 4, \label{eq:404a}\\
m\equiv & 10+8n  \pmod {25}, \label{eq:404b}\\
m\equiv &
 \begin{cases}
      1 \pmod 3 & \text{ if $n_3\not\equiv 2 \pmod 3$},  \\
      2 \pmod 3 & \text{ if $n_3\equiv 2  \pmod 3$}. \\
  \end{cases}
  \label{eq:404c}
\end{align}
Using condition \eqref{eq:404a}, we find that 
$$ 2\nmid L_1L_3(m), \ \ 2^1\|L_2(m).$$
Using condition \eqref{eq:404b}, we find that
$$ 5\nmid L_1L_3(m), \ \ 5^1||L_2(m).$$
Using condition \eqref{eq:404c}, we find that 
$$ 3\nmid L_1L_3(m), \ \ 3^1||L_2(m).$$

We deduce that  $L_1$ and $L_3$ have no fixed divisors  and that
$L_2$ has the fixed divisors $2,3,5$ with 
exponent one.
We adjoin extra factors of $11^2\cdot 17$ to $L_1$, $13^2\cdot 19$ to $L_3$, and then 
we apply the Basic Theorem to the reduced forms
\begin{equation*}
\frac{L_1(m)}{11^2 \cdot 17}, 
\frac{L_2(m)}{2\cdot 3 \cdot 5}.
\frac{L_3(m)}{13^2 \cdot 19},
\end{equation*}
Multiplying the original forms by $2$, $3$ and $5$ according to
\eqref{eq:401},  we arrive in all cases at numbers with 
exponent pattern $\{2,1,1,1,1\}$. 

{\bf Case 2.} For this case, we assume that Hypothesis $T$ holds for some 
prime $p>3$. We write $p=2k+1$, so that 
\begin{equation*}
p\equiv k \equiv 2 \pmod 3.
\end{equation*}
In this case, we will use the three linear forms
\begin{equation} \label{eq:410}
L_1(m)= 2m+n, \ \ L_2(m)=pm+kn, \ \ L_3(m)=(p+2)m+(k+1)n.
\end{equation}
Note the relations
\begin{equation} \label{eq:411}
pL_1-2L_2= (p+2)L_1-2L_3 = pL_3-(p+2)L_2 = n.
\end{equation}
We may assume that $3\cdot 5|n$; otherwise, we may appeal to Case 1.
We will always assume that 
\begin{equation}\label{eq:412}
m\equiv 1 \pmod 3;
\end{equation}
the motivation is that we ensure that $3\nmid L_1L_2L_3(m)$. 

{\bf Subcase 2a}. Assume that neither $p$ nor $p+2$ divides $n$. 
We also assume \eqref{eq:412}, and we specify that 
\begin{align} 
m \equiv  \ & 2+(k+1)n \pmod 4, \label{eq:412a} \\
L_3(m) \equiv \ & p \pmod {p^2}, \label{eq:412b} \\
m \equiv \ & 0 \pmod {p+2}. \label{eq:412c}
\end{align}

From \eqref{eq:412a}, we see that
\begin{align*} 
L_3(m) \equiv  & (p+2)m + (k+1)n \equiv 2p+ (p+3)(k+1)n \\
 \equiv & 2+2(k+2)(k+1)n \equiv 2 \pmod 4.
 \end{align*}
In the last line, we use the fact that $(k+1)(k+2)$ is always even. 
Furthermore, 
\begin{equation*}
2m+ n \equiv n \equiv 1 \pmod 2, \ \ pm+kn\equiv m+k \equiv (2k+1)n \equiv 1 \pmod 2.
\end{equation*}
Therefore, $2^1||L_3(m)$ and $2\nmid L_1L_2(m)$. 

From \eqref{eq:412b}, we see that $p^1|| L_3(m)$.  
Moreover, $p\nmid L_1L_2(m)$ because
\begin{align*}
L_1(m) = & 2m+n \equiv -kn \not\equiv 0 \pmod p, \text{ and } \\
L_2(m) = & pm+kn\equiv kn \not\equiv 0 \pmod p.
\end{align*}
From \eqref{eq:412c}, we see that $L_i(m)\equiv b_i n \pmod{p+2}$, where 
$b_1=1,b_2=k, b_3=k+1$. Therefore $(p+2)\nmid L_1L_2L_3(m)$.
We conclude that $L_3$ has a fixed divisor $2p$ and $L_1,L_2$ have no fixed divisors.

Let $p_1, p_2$ be distinct odd primes%
\footnote{We are no longer using the convention that $p_i$ denotes the $i^{\text{th}}$ prime.}
with $(p_1p_2, p(p+2)n)=1$. We adjoin extra factors of 
$p_1^2$ and $p_2^2$ to $L_1$ and $L_2$ respectively, then 
we apply the Basic Theorem to the reduced forms
\begin{equation*}
\frac{L_1(m)}{p_1^2}, \ \ 
\frac{L_2(m)}{p_2^2}, \ \
\frac{L_3(m)}{2p}.
\end{equation*}
Multiplying the original forms by $2, p, p+2$ according to \eqref{eq:411}, we obtain infinitely 
many positive integers $x$ such that $x$ and $x+n$ both have the exponent pattern $\{2,1,1,1\}.$
To get the desired exponent pattern $\{2,1,1,1,1\}$, we adjoin one more prime factor to each form.

{\bf Subcase 2b.} Assume that $p\nmid n$ and $(p+2)| n$, and write $n=(p+2)n_0$. 
We also assume \eqref{eq:412}, and we specify that 
\begin{align} 
m \equiv & \  2+3(p+2)(k+1)n \pmod 4,  \label{eq:413a}\\
m \equiv & \ -(k+1)n_0 + p \pmod{p^2},\label{eq:413b}\\
m(m+(k+1)n_0)  \not\equiv & \ 0 \pmod{p+2}.\label{eq:413c}
\end{align}
The last condition is possible because there are at most  two residue classes mod $(p+2)$ for
which \eqref{eq:413c} fails. 
Furthermore, from \eqref{eq:413a}, we see that
\begin{align*}
L_3(m) & =  (p+2)m+(k+1)n \equiv 2(p+2) + \left( 3(p+2)^2+1 \right)(k+1)n \\ &
    \equiv 2 \pmod 4.
\end{align*}
In the last equation, we use the fact that $j^2\equiv 1 \pmod 4$ whenever $j$ is odd.
We also note that
\begin{equation*}
L_1(m)=  2m+n \equiv    1   \pmod 2,  \ \ \text{and }
L_2(m)=  pm+kn \equiv   m+k  \equiv 1  \pmod 2 
\end{equation*}
Therefore $2^1||L_3(m)$ and $2\nmid L_1L_2(m)$.

From \eqref{eq:413b}, we see that 
\begin{equation*}
L_3(m) \equiv  (p+2) \left( m+(k+1)n_0 \right) \equiv  2p \pmod {p^2}.
\end{equation*}
Furthermore, 
\begin{align*}
L_1(m) \equiv -2kn_0 \pmod p \text{ and } 
L_2(m) \equiv kn  \pmod p.
\end{align*}
Therefore $p^1||L_3(m)$ and $p\nmid L_1L_2(m)$.
From \eqref{eq:413c}, we see that $(p+2)\nmid L_1L_2(m)$, and 
\begin{equation*}
(p+2)^1 || L_3(m) =(p+2)(m+(k+1)n_0).
\end{equation*}
We conclude that $L_3$ has fixed divisors $2p(p+2)$ and $L_1, L_2$ have no fixed
prime divisors. 

Let $p_1, p_2, p_3, p_4$ be distinct odd primes 
with $(p_1p_2p_3p_4, p(p+2)n)=1$. 
We adjoin extra factors of 
$p_1^2p_2$ and $p_3^2 p_4$ to $L_1$ and $L_2$ respectively, then 
we apply the Basic Theorem to the reduced forms
\begin{equation*}
\frac{L_1(m)}{p_1^2 p_2}, \ \ 
\frac{L_2(m)}{p_3^2 p_4}, \ \
\frac{L_3(m)}{2p(p+2)}.
\end{equation*}
Multiplying the original forms by $2, p, p+2$ according to \eqref{eq:411}, we obtain infinitely 
many positive integers $x$ such that $x$ and $x+n$ both have the exponent pattern $\{2,1,1,1,1\}.$


{\bf Subcase 2c.} Assume that $p| n$ and $(p+2)\nmid n$, and write $n=pn_0$. 
This case is similar to Subcase 2b except that the roles of $p$ and $p+2$ have been 
reversed. 
We  assume \eqref{eq:412}, and we specify that 
\begin{align} 
m \equiv & \ 2+3pkn \pmod 4,  \label{eq:414a}\\
m \equiv & \  -kn_0 + p+2  \pmod{(p+2)^2},\label{eq:414b}\\
m(m+k n_0)  \not\equiv & \ 0 \pmod{p}.\label{eq:414c}
\end{align}
From \eqref{eq:414a}, we see that
\begin{equation*}
L_2(m)  =  pm+kn \equiv 2p + \left( 3p^2+1 \right)kn  \equiv 2 \pmod 4.
\end{equation*}
Furthermore,
\begin{align*}
L_1(m)= & 2m+n \equiv 1   \pmod 2,  \text{ and } \\
L_3(m)= & (p+2)m+(k+1)n \equiv m+k +1 \equiv 1 \pmod 2.
\end{align*}
Therefore $2^1||L_2(m)$ and $2\nmid L_1L_3(m)$.

From \eqref{eq:414b}, we see that 
\begin{equation*}
L_2(m) \equiv  p \left( m+kn_0 \right) \equiv  p(p+2)  \pmod {(p+2)^2}.
\end{equation*}
Furthermore, 
\begin{align*}
L_1(m) \equiv & (-2k-2)n_0  \pmod {(p+2)},  \text{ and } \\
L_3(m) \equiv & (p-2k)n_0 \equiv n_0  \pmod {(p+2)}.
\end{align*}
Therefore $(p+2)^1||L_2(m)$ and $(p+2)\nmid L_1L_3(m)$.
From \eqref{eq:414c}, we see that 
\begin{equation*}
p \nmid L_1L_3(m) \text{ and }
p^1 || L_2(m) =p(m+kn_0).
\end{equation*}
We conclude that $L_2$ has fixed divisors $2p(p+2)$ and $L_1, L_3$ have no fixed
prime divisors. 

Let $p_1, p_2, p_3, p_4$ be distinct odd primes 
with $(p_1p_2p_3p_4, p(p+2)n)=1$. 
We adjoin extra factors of 
$p_1^2p_2$ and $p_3^2 p_4$ to $L_1$ and $L_3$ respectively, then 
we apply the Basic Theorem to the reduced forms
\begin{equation*}
\frac{L_1(m)}{p_1^2 p_2}, \ \ 
\frac{L_2(m)}{2p(p+2)},\ \ 
\frac{L_3(m)}{p_3^2 p_4}.
\end{equation*}
Multiplying the original forms by $2, p, p+2$ according to \eqref{eq:411}, we obtain infinitely 
many positive integers $x$ such that $x$ and $x+n$ both have the exponent pattern $\{2,1,1,1,1\}.$

\begin{remark} By considering Case 1a and Case 2a, we see that we also have the following  result.

\begin{theorem} \label{th:13a}
Let  $n$ be an odd number. Suppose that there exists a twin prime pair 
$(p,p+2)$ such that $p\nmid n$ and $(p+2)\nmid n$. 
Then there are infinitely many positive integers $x$ such that $x$ and $x+n$ both have
exponent pattern 
$$ \{2,1,1,1\}.$$
\end{theorem}

Note that this generalizes Theorem \ref{th:123}. Note also that using the set $\cT_1$ from 
\eqref{eq:T1}, we find that the hypothesis of Theorem \ref{th:13a} is true for any odd $n<S'$,
where 
$$ S'=\prod_{p\in \cT_1} p > 10^{10^{14}}.$$
\end{remark}


\section{The Case of Odd Shift--Proofs of Theorem \ref{th:12a}} \label{sec:5}

 We first do  the proofs for the results on $\omega$ 
and $d$; i.e., the proofs of \eqref{eq:12o} and \eqref{eq:12d}. 
The proof of \eqref{eq:12O} is much simpler, and it
will be done at the end of this section.
For the proofs of \eqref{eq:12o} and \eqref{eq:12d}, we will 
assume that $3\cdot 5 \cdot 7|n$; otherwise, the desired results 
 follow from Theorem \ref{th:11a}.

Consider the admissible system
\begin{equation} \label{eq:242}
672\ell+41, \ \  672\ell+47, \ \ 672\ell+55.
\end{equation}
By the Basic Theorem, there are infinitely many positive integers $\ell$ such that 
at least two of the forms in \eqref{eq:242} 
are $E_2$-numbers with both prime factors exceeding $C$. 
In this case, we take 
$C$ to be the greatest prime factor of $ 13n.$
In particular,
\begin{equation} \label{eq:C}
C\ge 13.
\end{equation}

We divide our argument into three cases.

{\bf Case 1:} $672\ell+41$ and $672\ell+47$ are both $E_2$-numbers. 
Let $q=672\ell+41$, and let $k=112\ell+7$, so that $q=6k-1$, $q\equiv 1\pmod 4$, and 
$k\equiv 3\pmod 4$.
Consider the system of linear forms
\begin{equation} \label{eq:243z}
L_1(m)= 6m+n, \ \ L_2(m)=qm+kn, \ \ L_3(m)=(q+6)m+(k+1)n.
\end{equation}
Note that 
\begin{equation} \label{eq:243}
6L_2-qL_1=6L_3-(q+6)L_1=(q+6)L_2-qL_3=n.
\end{equation}

Let $q_1, q_2$ be the prime divisors of $q$. We are assuming $3|n$, and we write 
$n=3n_3$. We specify that $m$ satisfies the conditions
\begin{align}
m \equiv & 2 \pmod 4, \label{eq:243a}\\
m(2m+n_3)  \not\equiv  & 0 \pmod 3, \label{eq:243b}\\
L_3(m) \equiv & q_1 \pmod {q_1^2}.\label{eq:243c}
\end{align}

From \eqref{eq:243a}, we find that 
\begin{equation*}
L_3(m)=(q+6) m + (k+1)n \equiv 3m \equiv 2 \pmod 4,
\end{equation*}
and that $2\nmid L_1L_2(m)$. 
From \eqref{eq:243b}, we find  that
\begin{equation*}
L_1(m) = 3(2m+n_3),
\end{equation*}
and so $3^1||L_1(m)$. Moreover, $3\nmid L_2L_3(m)$. 
From \eqref{eq:243c}, we find that $q_1^1||L_3(m)$.
Moreover,  from \eqref{eq:243} and \eqref{eq:243c}, we see that
\begin{equation*}
-6 L_1(m) \equiv 6L_2(m) \equiv n \pmod{q_1}.
\end{equation*}
Now $q_1$ exceeds all prime divisors of $n$, so $q_1\nmid L_1L_2(m)$.

To summarize, our congruence conditions \eqref{eq:243a} through \eqref{eq:243c} imply that
$L_1$ has fixed divisor $3$, 
$L_2$ has no fixed divisors,
and $L_3$ has fixed divisors $2$ and $q_1$. 

Let $p_1,p_2,p_3,p_4$ be distinct primes with
$$(p_1p_2p_3p_4,6q(q+6)n)=1.$$
We adjoin extra factors of  $p_1^2$ to $L_1$, $p_2^2p_3$ to $L_2$, and $p_4$ to $L_3$. 
We apply the Basic Theorem to the forms
\begin{equation*}
\frac{L_1}{3p_1^2}, \frac{L_2}{p_2^2p_3}, \frac{L_3}{2p_4q_1}.
\end{equation*}
We use the multipliers in \eqref{eq:243}, and we note that
\begin{align}
d(3p_1^2 q)= & d(3p_1^2(q+6)) = d(6p_2^2 p_3)= d((q+6)p_2^2p_3) \notag \\
= & d(6\cdot 2 q_1 p_4)=d(2q_1^2 q_2 p_4)=24. \label{eq:243d}
\end{align}
We  deduce that there are infinitely many $x$ such that
\begin{equation*}
d(x)=d(x+n)=96. 
\end{equation*}
By the usual procedure of adjoining  extra prime factors, we deduce that for any positive
integer $A$, $d(x)=d(x+n)=96A$ infinitely often. 
Moreover, we note that
\begin{align}
\omega(3p_1^2 q)= & \omega(3p_1^2(q+6)) = 
    \omega(6p_2^2 p_3)= \omega((q+6)p_2^2p_3) \notag \\
= & \omega(6\cdot 2 q_1 p_4)=\omega(2q_1^2 q_2 p_4)=4. \label{eq:243e}
\end{align}
We thereby obtain infinitely many $x$ such that $\omega(x)=\omega(x+n)=6$.

{\bf Case 2.} $672\ell+41$ and $672\ell+55$ are both $E_2$-numbers. 
This is similar to Case 1, so we will leave many of the details. 
The primary difference is that the prime $7$ plays the role played by $3$ in Case 1.

In this case, we let $k=48\ell+3$, so that $q=14k-1$ and $k\equiv 3 \pmod 4$.
Consider the system of linear forms
\begin{equation*}
L_1(m)=14m+n, \ \ L_2(m)=qm+kn, \ \ L_3(m)=(q+14)m+(k+1)n.
\end{equation*}
Note that
\begin{equation} \label{eq:245}
14L_2-qL_1= 14L_3-(q+14)L_1=(q+14)L_2-qL_3 = n.
\end{equation}

Let $q_1, q_2$ be the divisors of $q$. We are assuming that $7|n$, and we write $n=7n_7$. 
We specify that $m$ satisfies the conditions
\begin{align}
m\equiv \ & 2 \pmod 4, \label{eq:244a} \\
m(2m+n_7) \not\equiv \ & 0 \pmod 7, \label{eq:244b}\\
L_3(m) \equiv \ & q_1 \pmod{q_1^2}. \label{eq:244c}
\end{align}
From these conditions, we find that $2^1||L_3(m)$ and $2\nmid L_1L_2(m)$. 
We also find that $7^1||L_1(m)$ and $7\nmid L_2L_3(m)$. 
Moreover, $q_1^1||L_3(m)$ but $q_1\nmid L_1L_2(m)$. 
In summary, our congruence conditions imply that 
$L_1$ has fixed divisor $7$, $L_2$ has no fixed divisors, and $L_3$ has 
fixed divisors $2$ and $q_1$, similarly to Case 1. 

Let $p_1, p_2, p_3, p_4$ be distinct primes with 
$$ (p_1p_2p_3p_4,14q(q+6)n)=1.$$
We adjoin extra factors of $p_1^2$ to $L_1$, $p_2^2p_3$ to $L_2$ and $p_4$ to 
$L_3$. We apply the Basic Theorem to the forms
\begin{equation*}
\frac{L_1}{7p_1^2}, \ \ \frac{L_2}{p_2^2p_3}, \ \ \frac{L_3}{2p_4 q_1},
\end{equation*}
and we use the multipliers in \eqref{eq:245}. 
In this case, we note that \eqref{eq:243d} and \eqref{eq:243e} are true
when $6$ is replaced by $14$ and $3$ replaced by $7$. 
Accordingly, we deduce that for any positive integer $A$,
there are infinitely many $x$ such that
$ d(x)=d(x+n)=96A$ and $\omega(x)=\omega(x+n)=6$.


{\bf Case 3:} $672\ell+47$ and $672\ell+55$ are both $E_2$-numbers. 
Let $q=672\ell+47$, and let $k=84\ell+6$, so that $q=8k-1\equiv 3 \pmod 4$ and $k\equiv 2\pmod 4$.
Consider the system of linear forms
\begin{equation} \label{eq:247}
L_1(m)= 8m+n, \ \ L_2(m)=qm+kn, \ \  L_3(m)=(q+8)m+(k+1)n.
\end{equation}
Note that 
\begin{equation} \label{eq:244}
8L_2-qL_1=8L_3-(q+8)L_1=(q+8)L_2-qL_3=n.
\end{equation}

Let $q_1, q_2$ be the prime divisors of $q$, and choose $m$ so that 
\begin{align}
m\equiv & 4+(k+1)n \pmod 8, \label{eq:245a} \\
L_3(m) \equiv & q_1 \pmod {q_1^2},  \label{eq:245b} \\
L_3(m) \equiv & q_2^2 \pmod {q_2^3}. \label{eq:245c}
\end{align}

From \eqref{eq:245a}, we see that
\begin{equation*}
L_3(m)=(q+8)m+(k+1)n \equiv -m+(k+1)n \equiv 4 \pmod 8,
\end{equation*}
so $2^2||L_3(m)$.
Note also that $m$ and $n$ are odd and $k$ is even,  so $2\nmid L_1L_2(m).$ 
Conditions \eqref{eq:245b} and \eqref{eq:245c} imply that $L_3(m)$ has the fixed divisors $2^2,q_1,q_2^2$.  The forms $L_1$ and $L_2$ have no fixed divisors, similarly to 
Cases 1 and 2. 

Let $p_1,p_2,p_3,p_4$ be primes with 
\begin{equation} \label{eq:244k}
(p_1p_2p_3p_4, 2q(q+8)n)=1.
\end{equation}
We add extra divisors $p_1^2 p_2^2$ to $L_1$ and $p_3^2 p_4^2$ to $L_2$. 
We apply the Basic Theorem to the forms
\begin{equation*}
\frac{L_1}{p_1^2p_2^2}, \frac{L_2}{p_3^2 p_4^2}, \frac{L_3}{2^2 q_1 q_2^2}.
\end{equation*}
We use the multipliers in \eqref{eq:244} and note that 
\begin{align*}
d(p_1^2 p_2^2 q) & = d(p_1^2 p_2^2 (q+8)) = d(8p_3^2 p_4^2) = d(p_3^2 p_4^2 (q+8)) \\
 & = d(2^3\cdot 2^2 q_1 q_2^2) = d(2^2 q_1 q_2^2 q) = 36.
 \end{align*}
 We  deduce that there are infinitely many $x$ with
 \begin{equation*}
 d(x)=d(x+n)=144.
 \end{equation*}
By adjoining further prime factors, we see that for any positive
integer $A$, $d(x)=d(x+n)=144A$ infinitely often. 

The above argument has to be modified slightly for $\omega$. We take the same linear forms
as given in \eqref{eq:247}. As before, we take $m\equiv 4+(k+1)n \pmod 8$, and we let $q_1, q_2$ be the 
prime divisors of $q$. Let $q_3, q_4$ be the prime divisors of $q+8$. We choose $m$ so that
\begin{align*}
L_2(m) \equiv q_3 \pmod {q_3^2}, \text{ and }
L_3(m) \equiv q_1 \pmod {q_1^2}.
\end{align*}
Let $p_1,p_2$ be as in \eqref{eq:244k}. 
Adjoin an extra factor of $p_1$ to $L_1$ and an extra factor of $p_2$ to $L_2$, then 
apply the Basic Theorem to the forms
\begin{equation*}
\frac{L_1}{p_1}, \frac{L_2}{p_2q_3}, \frac{L_3}{2^2q_1q_2}.
\end{equation*}
We use the multipliers in \eqref{eq:244} and note that 
\begin{align*}
\omega(p_1 q)& = \omega(p_1(q+8)) = \omega(8p_2q_3) = \omega(p_2 q_3 (q+8))\\
& =  \omega(8\cdot 2^2 q_1 q_2) = \omega(2^2 q_1 q_2 q)=3.
\end{align*}
We deduce that there are infinitely many $x$ such that 
\begin{equation*}
\omega(x)=\omega(x+n)=5.
\end{equation*}
To make this compatible with Cases 1 and 2, we adjoin one further prime factor
so that we obtain infinitely many $x$ with $\omega(x)=\omega(x+n)=6$. 
This concludes Case 3.

We can now prove \eqref{eq:12o} and \eqref{eq:12d}. From the above three cases, we
we see that there are infinitely many $x$ such that $\omega(x)=\omega(x+n)=6$. 
Statement \eqref{eq:12o} follows from the usual procedure of adjoining further prime factors.
For \eqref{eq:12d}, we note that Cases 1 and 2 give infinitely many $x$ such that 
$d(x)=d(x+n)=96A$, while Case 3 gives $d(x)=d(x+n)=144A$. Statement \eqref{eq:12d}
follows upon noting that $\text{lcm}(96,144)=288$.

For the corresponding result on $\Omega$, we return to the simpler procedure used
in Section \ref{sec:3}. 
First, we note that by Theorem \ref{th:13a}, we have \eqref{eq:12O} with $A=5$ whenever
$(n,15)=1, 3, $ or $5$. Henceforth, we assume that $15|n$.

We take the forms given in \eqref{eq:141}; i.e., 
\begin{equation*}
L_1(m)= 2m+n, \ \ L_2(m)=3m+n, \ \ L_3(m)=5m+2n.
\end{equation*}
We  may write $n=3n_3=5n_5$. We specify that
\begin{align*}
m\equiv  \  & n+2 \pmod 4, \\
m(m+n_3) \not\equiv\ & 0 \pmod 3,\\
m(m+2n_5) \not\equiv\  & 0 \pmod 5.
\end{align*}
From these conditions, we see that $L_1$ has no fixed divisors, $L_2$ has fixed divisors $2$ and $3$, and $L_3$ has fixed divisor $5$. We apply the Basic Theorem to the reduced forms
\begin{align*}
L_1(m), \ \ \frac{L_2(m)}{2\cdot 3}, \ \ \frac{L_3(m)}{5}.
\end{align*}
We adjoin one new prime factor to $L_3$ and two new prime factors to $L_1$. 
We then find that there are at least two forms $L_i$ such that  
$$\Omega(L_i(r))=4$$
for infinitely many $r$. Using the relations in \eqref{eq:141a}, we find that there
are infinitely many $x$ such that
$$\Omega(x)=\Omega(x+n)=5.$$

The completes the proof of \eqref{eq:12O} in the case $A=5$. 
The general case follows by adjoining further prime factors.

\footnotesize
D. A. Goldston,
Department of Mathematics,
San Jose State University,
San Jose, CA 95192, USA,
e-mail: goldston@math.sjsu.edu

S. W. Graham,
Department of Mathematics,
Central Michigan University,
Mt. Pleasant, MI 48859, USA,
email: sidney.w.graham@cmich.edu

J. Pintz,
R\'enyi Mathematical Institute of the Hungarian Academy of Sciences,
H-1053 Budapest, Realtanoda u. 13--15.
Hungary,
e-mail: pintz@renyi.hu

C. Y. Y{\i}ld{\i}r{\i}m, Department of Mathematics,
Bo\~{g}azi\c{c}i University, Bebek, Istanbul 34342, Turkey, \&
Feza G\"ursey Enstit\"us\"u \c Cengelk\"oy, Istanbul, P.K. 6,
81220 Turkey,\\
e-mail: yalciny@boun.edu.tr

\end{document}